\documentclass[12pt]{article}

\usepackage{amssymb,latexsym}

\title{ On Number of Compositions of Natural Numbers}
\author{\bf Milan Janji\'c\\Department of Mathematics and Informatics,\\ University of Banja Luka, Republic of Srpska}

\begin{document}
\maketitle

\begin{abstract}
 We first give a combinatorial interpretation of coefficients of Chebyshev polynomials, which allows to connect them with compositions of natural numbers.
    Then we describe a relationship between the number of compositions of a natural number in which a certain number of  parts are   $p-1,$ and other parts are $\geq p,$ with compositions in which all parts are $\geq p.$
 Then we find a relationship between  principal minors of a type of Hessenberg matrices and compositions of natural numbers.
\end{abstract}

\section{Introduction}

  We first  prove some combinatorial results which will be used in the paper. In the second section we prove a formula for some restricted compositions of natural numbers.

  In the second section  we show that coefficients of Chebyshev  polynomials  may be interpreted in terms of compositions of natural  umbers.

   Next, in the third section,  we connect compositions in which some parts are $p-1$, and all other parts $\geq p$ with compositions in which all parts are $\geq p.$

Different mathematical objects may be represented by determinants of  Hessenberg matrices. We describe a special kind of Hessenberg matrices which determinants are compositions of natural numbers. It will also be shown that sums of principal minors of these matrices, that is, coefficients of its characteristic polynomials also represent some type of compositions of natural numbers. This is done in the fourth section.

 We shall use the following result about
Hessenberg's matrices, which is easily proved by induction.

\noindent\textbf{Theorem 1.} \textit{ Let $a_1,p_{i,j},\;( i\leq j)$ be arbitrary elements of a commutative ring $R,$ and let the sequence $a_1,a_2,\ldots$ be defined by
\[a_{n+1}=\sum_{i=1}^np_{i,n}a_i,\;(n=1,2,\ldots),\]
If
\[P_n=
\left[\begin{array}{llllll}
p_{1,1}&p_{1,2}&p_{1,3}&\cdots&p_{1,n-1}&p_{1,n}\\
-1&p_{2,2}&p_{2,3}&\cdots&p_{2,n-1}&p_{2,n}\\
0&-1&p_{3,3}&\cdots&p_{3,n-1}&p_{3,n}\\
\vdots&\vdots&\vdots&\ddots&\vdots&\vdots\\
0&0&0&\cdots&p_{n-1,n-1}&p_{n-1,n}\\
0&0&0&\cdots&-1&p_{n,n}\end{array}\right],\] then
\[a_{n+1}=a_1\det A_n,\;(n=1,2,\ldots).\]}

  We prove now  a combinatorial result which will also  be used in the paper.

   Consider a set $X$ consisting of $n$ blocks $X_i$ each having $2$ elements,
 and an additional block $Y$ with $m$ elements.  We call $X_i$ the main blocks of $X.$  By an  $n+k$-inset of $X$  we shall
 mean a subset of $X$ with $n+k$ elements intersecting each main block $X_i,\;(i=1,2,\ldots,n).$ We shall denote by $N(n,k,m)$ the number of $n+k$-insets of $X.$
For this number  we have

\noindent\textbf{Proposition 1.} \textit{Let $n$ be a positive integer, and let $k,m$ be nonnegative integers. Then
 \[N(n,k,m)=\sum_{i=0}^n(-1)^i{n\choose i}{2n-2i+m\choose n+k}.\]}
Also,
\[N(n,k,m)=2^{n-k}\sum_{i=0}^m2^{i}{m\choose i}{n\choose k-i}.\]

{\it Proof.}
For $i=1,2,\ldots,n$
 and an $n+k$-subset $Z$ of $X$ define the property $i$ to be:

\begin{center} $Z$ does not intersect $X_i.$ \end{center}

By inclusion-exclusion principle we obtain \[f(n,k,m)=\sum_{I\subseteq [n]}(-1)^{|I|}N(I),\] where
$N(I)$ is the number of $n+k$-subsets of $X$ which do not intersect main blocks $X_i,\;(i\in I).$ There are
\[{n\choose i}{2n+m-2i\choose n+k}\] such subsets and the first  formula is proved.

We may count $n+k$-insets of $X$ in another way.
We prove firstly the following formula:
 \begin{equation}\label{ff}N(n,k,m)=\sum_{i=0}^m{m\choose i}N(n,k-i,0).\end{equation}

In fact, omitting the additional block  $Y$ we
obtain the set $X_1$ with no additional block.

Each $n+k$-inset of $X_1$ is an $n+k$-inset of $X.$  There are $N(n,k,0)$ such insets. In this way we obtain all
$n+k$-inset of $X$ not intersecting $Y.$

The remaining $n+k$-insets of $X$ are obtained as a union of  some $n+k-i,\;(1\leq i\leq m)$-inset of $X_1$ and some
$i$-set of additional block $Y$. This proves (\ref{ff}).

Further, we obtain $N(n,k,0)$ by choosing two elements from each of arbitrary $k$ main blocks,
 and one element from each of remaining $n-k$ main blocks. The first may be done in $n\choose k$ ways and the second in $2^{n-k}$ ways.
We conclude that $N(n,k,0)=2^{n-k}{n\choose k}.$ It follows that

\[N(n,k,m)=2^{n-k}\sum_{i=0}^m2^{i}{m\choose i}{n\choose k-i},\] and the second formula is true.\hfill$\Box$

\section{Chebyshev Polynomials and Compositions}
In this section we shall show that  coefficients of Chebyshev polynomials
may be interpreted in terms of compositions of natural numbers.

Let $x$ be a symbol. By an  $u(n,k,x)$-sequence we shall mean a sequence of length $n$ in which $k$ terms are equal $x,$ and remaining terms are either $1$ or $1|.$ The sequence does not end by $1|.$

\noindent\textbf{Proposition 2.} \textit{The number  of $u(n,k,x)$-  sequences is $N(n-1,k,1).$}

{\it Proof.}
 We prove  that all  $u(n,k,x)$-  sequences may be labeled by $n-1+k$-insets of a set $X$ consisting of $n-1$ main blocks, each with $2$ elements, and an additional block $Y$  with one element.

Let $Z$ be an $n-1+k$-inset of $X.$ Form a $u(n,k,x)$ sequence as follows. Put either  $x$ at the last place, if the element from $Y$   is in $Z,$ or put $1$ if the element from $Y$ is  not in $Z.$ Put then the remaining $x$'s on
positions corresponding to blocks from which both elements are in $Z.$ Then insert $1$'s on the remaining places to obtain a $n$-sequence  in which  $k$ elements are $x$'s, and the remaining $n-k$ elements are $1$'s.
The last element in the sequence  is either $x$ or $1.$

If the last elements is $x$ then we have $n-k$ main blocks from which only one element is in $Z.$
 Put these blocks and all $1'$s into a  bijective correspondence.
Each   $1$, is either followed by $|,$ if the second element of corresponding block is in $Z$,  or not followed by $|,$ if the first element of corresponding block is in $Z.$

If the sequence ends by $1$ then there are $n-1-k$ main blocks from which one element is in $Z.$ We now put these block into a bijective correspondence with $1$'s, excluding the last one. In such way we also obtain a $u(n,k,x)$-sequence.
It follows that each $n-1+k$-inset of $X$ determines a $u(n,k,x)$-sequence.
Clearly,  that different insets produce different sequences.

Conversely, consider a $u(n,k,x)$-sequence. Form a $n-1+k$ inset $Z$  of $X$ in the following way: If the last element in the sequence is $x$ then put element of $Y$ in $Z.$ Otherwise this element is not in $Z.$

Furthermore, put in $Z$  both elements of main blocks which indices correspond to the indices of remaining  $x$ in $u(n,k,x).$
 The indices of the remaining main blocks correspond to $1$'s, excluding eventually the last $1.$ If $1$ corresponding to such main block has $|$ then put the second element of the block in $Z.$   Otherwise put the first element in $Z.$ Thus, $Z$ is a $n-1+k$- inset of $X.$
\hfill$\Box$

In the paper [2] a family of polynomials $P_n(x)$ is defined such that the coefficient by $x^{n-k+m}$ of $P_{n+k+m}(x)$ is
 $(-1)^{k}N(n,k,m).$ Particulary, for $m=0$ we have Chebyshev polynomials  $U(n,x)$of the second kind, and for $m=1$ we obtain Chebyshev polynomials $T_{n,k}(x),$  of the first kind. Using Propositions 1 and 2 we obtain

\noindent\textbf{Proposition 3.} \textit{ Let $t_{n,k}$ be  the coefficient of Chebyshev polynomial  $T_{n+k+1}(x)$ by $x^{n-k+1}.$ Then $(-1)^kt_{n,k}$ is equal to the number of  $u(n+1,k,x)$- sequences.\\
Particulary, we have the following two explicit formulas for $t_{n,k}.$
\[(-1)^kt_{n,k}=\sum_{i=0}^n(-1)^i{n\choose i}{2n-2i+1\choose n+k}.\]
Also,
\[(-1)^kt_{n,k}=2^{n-k}\left[{n\choose k}+2{n\choose k-1} \right].\]}

There is a natural relationship between $u(n,k,x)$- sequences  and particular kind of compositions of natural numbers.
 We shall give one such relationship.
 In fact, if $a$ is a nonnegative integer then sequences $s(n,k,a|)$ produce compositions of $n+k(a-1),$ obtained by adding numbers not separated by $|.$ But, we want that different sequences give different compositions. This is not the case if $a\leq n-k,$
 since, for example, $11|2|$ and $2|11$ generate the same composition $2,2$ of $4.$

\noindent\textbf{Proposition 4.} \textit{Let $k,n$ be nonnegative integers such that $k\leq n+1.$ Then  $(-1)^kt_{n,k}$  is the number of compositions of $(n+1)(k+1)-k^2$ in which exactly $k$ parts are $\geq n-k+1.$ }

{\it Proof.} Denote $a=n-k+1.$ It is clear that each $s(n,k,a|)$-sequence produces a required composition. Different sequences give different
composition. Namely, consider the composition $(i_1,\ldots,i_m)$ obtained from two $u(n,k,a|)$- sequences.
Each  $i_t,\;(t=1,\ldots,m)$ may be obtained in two ways:
either from a term of the form $11\cdots 1|$ or from a term of the form  $11\cdots 1a|.$ Since no two such terms are equal it follows that $i_t$ must be obtain in the same way in both sequences. Hence, the sequences are the same.

Conversely, let $(i_1,i_2,\ldots,i_m)$ be a composition of $(n+1)(k+1)-k^2$ in which  exactly $k$ parts are $\geq n-k+1.$
Replace each such part by $11\cdots x,$ and all other parts by $11\cdots1|$ to obtain a $u(n,k,x)$ sequence.
The formula now follows from Proposition3.

In the case $k=1,\;n=3$ we have $a=3$ and the following $8$  compositions of $5$:
\[(3,2),(3,1,1),(4,1),(1,3,1),(5),(1,4),(1,1,3),(2,3).\]

\section{Some Restricted Compositions}

Let $n,p$ be  positive integers. We shall denote by $c(n,p)$ the number of compositions of $n$ in which all parts are $\geq p.$
It is convenient to define $c(0,p)=1.$
We know that  $c(n,1)=2^{n-1},\;c(n,2)=f_{n-1},$ where $f_{n-1}$ is Fibonacci number.

If, additionally, $k$ is a nonnegative integer then by $c(n,k,p,p-1)$ will be denoted the number of compositions of $n$ in which
  exactly $k$ parts are equal $p-1,$ and all other parts are $\geq p.$

The following result will be  a formula which connects $c(n+kp+1,k,p,p-1)$ with $c(n,p).$

\noindent\textbf{Proposition 4.} \textit{Let $n,p$ be  positive integers, and let $k$ be a nonnegative integer.  Then \begin{equation}\label{cnk}c(n+kp+1,k,p,p-1)=\sum_{j_1+j_2+\cdots+j_{k+1}=n}c(j_1+1,p)c(j_2+1,p)\cdots c(j_{k+1}+1,p),\end{equation}
where the sum is taken over $j_t\geq -1,\;(t=1,2,\ldots,k+1).$}

{\it Proof.}
We use the induction with respect of $k.$
For $k=0,$ the assertion is obviously true.

Assume that the assertion is true for $k-1.$ The greatest  value of $j_{k+1}$ is $n+k,$ and is obtained
for $j_1=j_2=\cdots=j_k=-1.$
Hence, we may write $(\ref{cnk})$ in the form:
\[c(n+nk+1,k,p,p-1)=\sum_{j=-1}^{n+k}c(j+1,p)\sum_{j_1+j_2+\cdots+j_{k}=n-j}c(j_1+1,p)\cdots c(j_{k}+1,p).\]
By the induction hypothesis we have
\begin{equation}\label{cnk1}c(n+kp+1,k,p,p-1)=\sum_{j=-1}^{n+k}c(j+1,p)\cdot c(n+kp+1-j-p,k-1,p,p-1).\end{equation}
Let $(i_1,i_2,\ldots)$ be a composition of $n$ with exactly $k$ parts equal $p-1$, and all other parts are $\geq p.$  Assume that the first $p-1$ is $i_m$.
Then $(i_1,i_2,\ldots,i_{m-1})$ is a composition of $j+1=i_1+\cdots+i_{m-1}$ in which all parts are $\geq p.$ Also,  $(i_{m+1},\ldots)$ is a composition of $n+kp+1-j-p$ with exactly $k-1$ parts equal $p-1,$ and all other parts $\geq p.$

It follows that the number of compositions in which the first $p-1$ is at the $m$th place is  $c(j+1,p)\cdot c(n-j-p,k-1,p,p-1).$ This is a term in the sum on the right side of (\ref{cnk1}).

For $m=1$ we have $j=-1$ which produces the first term in (\ref{cnk1}). If $p-1$  are on the last $k$ places of a composition  then
$j=n+k$ that gives  the last term in the sum of (\ref{cnk1}). Summing over all $j$ we conclude that the proposition is true.
 \hfill{$\Box$}

As an immediate consequence of Proposition 4, for $p=1$, we have

\noindent\textbf{Corollary 2.} \textit{ Let $n$ be a positive integer, and let $k$ be a nonnegative integer.
Then the number $K$ of weak compositions of $n+1$ in which exactly   $k$ parts are equal  $0$
is \[K=\sum_{j_1+j_2+\cdots+j_{k+1}=n}c(j_1+1,1)c(j_2+1,1)\cdots c(j_{k}+1,1),\] where $j_t\geq -1,\;(t=1,2,\ldots,k+1),$
and $c(0,1)=1,\;c(j,1)=2^{j-1},\;(j>0)$.}

Also, for $p=2$ we have

\noindent\textbf{Corollary 3.} \textit{ Let $n$ be a positive integer, and let $k$ be a nonnegative integer.
Then the number $K_1$ of composition of $n+2k+1$ in which exactly $k$ parts are equal  $1,$ and all other part are $\geq 2$
is \[K_1=\sum_{j_1+j_2+\cdots+j_{k+1}=n}f_{j_1}\cdots f_{j_{k+1}},\] where $j_t\geq -1,\;(t=1,2,\ldots,k+1),$ and $f_s$ are Fibonacci numbers.}

\section{Principal Minors of Some Hessenberg Matrices}
For positive integers $p,n$ we define an upper Hessenberg matrix
$F_{n,p}=(F_{n,p}(i,j))_{n\times n}$ in the following way:
\[F_{n,p}(i,j)=\left\{\begin{array}{rr}
-1&i=j+1,\\1&j-i\geq p-1\\0&\mbox{ otherwise }
\end{array}\right..\]

By Theorem 1 we easily obtain that
\[\det F_{1,p}=\det F_{2,p}=\det F_{p-1,p}=0,\;\det F_{p,p}=1,\]
and,
 \[\det F_{n,p}=\sum_{i=1}^{n-p+1}\det F_{i,p},\]
 that is,
\[\det F_{n,p}=\det F_{n-p,p}+\det F_{n-1,p}.\]
According to [1, p.63] we have

\noindent\textbf{Proposition 5.} \textit{ Let $n,p$ be positive  integers.
Then $c(n,p)=\det F_{n,p}.$}

Now, we shall find a relationship between coefficients of characteristic polynomials of matrices $F_{n,p}$ and  compositions mentioned in the preceding section.

\noindent\textbf{Proposition 6.} \textit{
Let $n,p$ be  positive integers, and let $k$ be a nonnegative integer. Then the number  $c(n+kp-2k,k,p,p-1)$  is the sum of all principal minors of the order $n-k$ of the matrix $F_{n,p}.$
}

{\it Proof.}
Denote  $M(i_1,\ldots,i_k)$ the minor of order $n-k$  obtained by deleting rows and columns of $F_{n,p}$ which indices are $1\leq i_1<i_2<\cdots<i_k\leq n.$
It is easy to see that the following equation holds
\[M(i_1,\ldots,i_k)=\det F_{i_1-1,p}\cdot\det F_{i_2-i_1-1,p}\cdots \det F_{i_{k}-i_{k-1}-1,p}\cdot\det F_{n-i_k,p}.\]

 According to Propositions 5  we have
\begin{equation}\label{minor}M(i_1,\ldots,i_k)= c(i_1-1,p)\cdot c(i_2-i_1-1,p)\cdots c(n-i_{k},p),\end{equation}

Denote $i_1-1=j_1+1,i_2-i_1-1=j_2+1,\ldots, i_k-i_{k-1}-1=j_k+1,\;n-i_{k}=j_{k+1}+1$ to obtain
that $j_1+\cdots+j_{k+1}=n-2k-1,\;(j_t\geq -1,\;(t=1,2,\ldots,k+1)).$

 It follows that

\[M(i_1,\ldots,i_k)= c(j_1+1,p)\cdot c(j_2+1,p)\cdots c(j_{k+1}+1,p).\]

Summing over all $1\leq i_1<i_2<\cdots<i_k\leq n$ we obtain a sum $S$ such that
\[S=\sum_{j_1+\cdots+j_{k+1}=n-2k-1} c(j_1+1,p)\cdot c(j_2+1,p)\cdots c(j_{k+1}+1,p),\]
where the sum is taken over all $j_t\geq -1,\;(t=1,2,\ldots,k+1).$

The proposition follows from Proposition 4.
\hfill$\Box$

 Since coefficients of the characteristic polynomial of a matrix are up to sign equal to the sum of its principal minors we may state Proposition 5 in the following form:

\noindent\textbf{Corollary 4.} \textit{If  $(-1)^{n-k}a_{n-k}$ is the coefficient by $x^k$ of the characteristic polynomial of the matrix  $F_{n,p}$
then $a_{n-k}$ is he number of compositions of $n+kp-2k,$ in which there are  $k$ parts equal $p-1$ and all other parts are $\geq p.$}

We state two particular cases of  Proposition 5, for $p=1$ and $p=2.$

\noindent\textbf{Corollary 5.} \textit{The number of weak compositions of $n-k,\;(k\leq n)$ in which exactly   $k$ parts are $0$, is equal to the sum of all principal minors of order $n-k$ of the matrix $F_{n,1}.$}

\noindent\textbf{Corollary 6.} \textit{The number of compositions of $n,$ in which exactly   $k,\;(k\leq n-2)$ parts are $1,$ and all other parts are $\geq 2,$  is equal to the sum of all principal minors of order $n-k$ of the matrix $F_{n,2}.$}

\noindent\textbf{References}

\vspace{0.2cm}

\noindent [1] G. E. Andrews {\it The Theory of Partitions,
} Addison-Wesley, 1976.

\noindent [2] M. Janji\'c, On a Class of Polynomials with Integer Coefficients, {\it Journal of Integer Sequences,} Volume 11, Issue 5, a:08.5.2, 2008.

\end{document}